\documentclass[onefignum,onetabnum]{siamonline250211}

\usepackage{lipsum}
\usepackage{amsfonts}
\usepackage{graphicx}
\usepackage{epstopdf}
\usepackage{algorithmic}
\usepackage{mathrsfs}

\ifpdf
  \DeclareGraphicsExtensions{.eps,.pdf,.png,.jpg}
\else
  \DeclareGraphicsExtensions{.eps}
\fi

\usepackage{enumitem}
\setlist[enumerate]{leftmargin=.5in}
\setlist[itemize]{leftmargin=.5in}

\newsiamremark{remark}{Remark}
\newsiamremark{hypothesis}{Hypothesis}
\crefname{hypothesis}{Hypothesis}{Hypotheses}
\newsiamthm{claim}{Claim}
\newsiamremark{fact}{Fact}
\crefname{fact}{Fact}{Facts}

\headers{Interneuron Subtypes and Bump Dynamics}{B. Ahmed, H. Cihak, and G. Handy}

\title{Stability and Wandering of Bumps in Neural Fields with Interneuron Subtypes\thanks{Submitted to the editors September 11, 2026.
\funding{This work was supported by the Burroughs Wellcome Fund’s Career Award at the
Scientific Interface to G. Handy.}}}

\author{Bilal Ahmed\thanks{School of Mathematics, University of Minnesota, Minneapolis, MN 
  (\email{ahme0692@umn.edu}, \email{cihak019@umn.edu},
  \email{ghandy@umn.edu}).}
\and Heather Cihak\footnotemark[2] \thanks{These authors contributed equally to this work}
\and Gregory Handy\footnotemark[2] \footnotemark[3]}

\usepackage{amsopn}

\ifpdf
\hypersetup{
  pdftitle={Stability and Wandering of Bumps in Neural Fields with Interneuron Subtypes},
  pdfauthor={B. Ahmed, H. Cihak, and G. Handy}
}
\fi

\begin{document}

\maketitle

\begin{abstract}
The maintenance of continuous variable information in working memory is thought to rely on persistent patterns of cortical activity. In delayed-estimation tasks, neural activity can form localized activity peaks, or ``bumps,'' whose positions track the remembered variable. Such activity is well described by continuous-attractor neural field models, but most existing models collapse cortical inhibition into a single homogeneous population. Here, we introduce a stochastic neural field model with distinct excitatory, parvalbumin-expressing (PV), and somatostatin-expressing (SST) populations to examine how inhibitory subtype structure shapes persistent activity. Using a Heaviside firing-rate approximation, we derive stationary bump solutions and reduce their linear stability to separate shifting and scaling modes. We show that population thresholds and inhibitory timescales determine both bump stability and the mechanism by which stability is lost, while inhibitory connection strengths and spatial scales substantially reshape the stable parameter region. In particular, broader SST connectivity promotes stable bump states. Finally, we derive an effective diffusion coefficient for noise-driven bump wandering and show that increasing the SST spatial footprint reduces the rate of memory diffusion. Together, these results demonstrate how inhibitory subtype structure can shape both the deterministic stability and stochastic precision of continuous-attractor memories.
\end{abstract}

\begin{keywords}
neural fields, working memory, interneuron subtypes, wandering, bump attractor
\end{keywords}

\begin{MSCcodes}
35B36, 35R60, 92B20
\end{MSCcodes}

\section{Introduction}
\label{sec:intro}
The ability to encode and maintain information for a few seconds is a defining feature of working memory, a core cognitive function that supports the planning and execution of everyday behavior~\cite{goldman1995cellular,constantinidis2016,curtis2006prefrontal}. Often conceptualized as a mental sketchpad, this transient system provides the spatial and temporal continuity between past experience and present action~\cite{fuster1973unit,goldman1995cellular, curtis2006prefrontal}. Delayed-response paradigms, such as the oculomotor delayed response task, are used investigate the phenomenon of visuo-spatial working memory and the link between neural mechanisms and information encoding~\cite{funahashi1989mnemonic, goldman1995cellular,ploner1999errors}. In such tasks, a subject briefly views a cue, encodes a continuous feature (e.g., position or orientation), maintains it through a delay period in the absence of the stimulus, and then reports it with a motor response~\cite{funahashi1989mnemonic, ploner1999errors, constantinidis2016}. Response errors accumulate during the delay, with error variance growing approximately linearly with delay length~\cite{white1994, ploner1998temporal,wimmer2014bump}. This time-dependent degradation motivates treating continuous-variable working memory as a dynamic, stochastic process~\cite{compte2000synaptic, kilpatrick2013wandering, bays2014, schneegans2018, panichello2019}.

A widely studied cellular component for this memory is persistent neural activity in the prefrontal cortex~\cite{fuster1973unit, funahashi1989mnemonic, wang1999synaptic, zylberberg2017, constantinidis2018persistent}. When a cue is remembered, specific populations of neurons sustain elevated firing throughout the delay period, even after the stimulus is removed~\cite{fuster1973unit, funahashi1989mnemonic, compte2000synaptic}. Neurons in the dorsolateral prefrontal cortex exhibit topographically organized ``memory fields,'' firing tonically only when a particular location is held~\cite{funahashi1989mnemonic, curtis2006prefrontal}, and this tuning is shared by neighboring putative interneurons, indicating that inhibition participates directly in shaping the mnemonic representation~\cite{rao1999, constantinidis2002}. At the population level, activity forms a localized peak, or ``bump,'' whose position tracks the remembered location~\cite{compte2000synaptic, wimmer2014bump}. The bump is sustained by recurrent excitation among pyramidal neurons and stabilized by feedback inhibition~\cite{wang1999synaptic, compte2000synaptic, wang2004}. Because the underlying network is marginally stable along the feature dimension, fluctuations arising from synaptic and spiking variability~\cite{faisal2008} cause the bump to drift stochastically~\cite{compte2000synaptic, kilpatrick2013wandering, constantinidis2018persistent, li2021drifts}. The variance of the bump position grows linearly over short delays~\cite{kilpatrick2013wandering, wimmer2014bump, constantinidis2018persistent}, providing a mechanistic account of the delay-dependent behavioral errors described above~\cite{white1994,ploner1998temporal, compte2000synaptic, wimmer2014bump, schneegans2018}.

Continuous neural field theory provides a tractable mathematical bridge between these biological observations and formal analysis~\cite{wilson1973,ermentrout1998, bressloff2012, coombes2014}. Notably,~\cite{amari1977dynamics} showed that a homogeneous field with local recurrent excitation and broader lateral inhibition supports localized activity bumps with exactly the marginal stability needed to encode a continuous variable anywhere along the feature axis. In the high-gain limit, where the firing rate becomes a Heaviside step, a bump is characterized entirely by its threshold-crossing interfaces, making existence and stability accessible in closed form~\cite{amari1977dynamics, coombes2004, coombes2012, faye2018}. This tractability carries over to the stochastic setting: incorporating noise and projecting onto the bump's neutral translational mode yields an effective diffusive description of wandering~\cite{bressloff2010, kilpatrick2013wandering, bressloff2015} that links network architecture directly to behavioral variability~\cite{kilpatrick2013interareal, kilpatrick2013heterogeneity}. The same framework has served as a testbed for many hypothesized memory mechanisms, including the encoding of certainty~\cite{carroll2014encoding, cihak2024robustly}, synaptic plasticity~\cite{seeholzer2019, eissa2022learning}, and serial dependence across trials~\cite{kiyonaga2017, kilpatrick2018, barbosa2020}, as part of a broader effort toward ``information-rich'' neural field models~\cite{zhou2025probabilistic} that relax strict assumptions while retaining mathematical tractability.

Across this body of work, inhibition is treated as a single homogeneous population, whether folded into a lateral-inhibition kernel acting on one excitatory layer~\cite{amari1977dynamics} or given its own dynamical layer in two-population excitatory-inhibitory (E/I) models~\cite{pinto2001b, laing2001, blomquist2005localized, folias2011}. Even within this simplified description, the properties of the inhibitory population turn out to be decisive for memory dynamics. In the E/I neural field of~\cite{cihak2022distinct}, the inhibitory firing threshold, timescale, and the strength and spatial extent of E/I coupling each set the boundaries of the stable bump regime as well as the type of instability encountered at those boundaries; bump position variance depends non-monotonically on these same parameters; and both stability and wandering are altered by inhibitory-to-inhibitory connections weak enough to be routinely neglected. Given this sensitivity to how a single inhibitory population is parameterized, it is natural to ask what is missed by ignoring the heterogeneity within that population.

Indeed, cortical interneurons do not form a single class; they comprise several molecularly, morphologically, and physiologically distinct GABAergic subtypes~\cite{rudy2011, abbas2018, tremblay2016, riedemann2019diversity, campagnola2022local}. While there are many interneuron subtypes, we focus on the two largest: parvalbumin-expressing (PV) and somatostatin-expressing (SST) interneurons~\cite{rudy2011, tremblay2016, riedemann2019diversity, campagnola2022local}. Importantly, these populations have been shown to play different structural and computational roles. Namely, PV cells provide fast, perisomatic inhibition by targeting the cell bodies of principal excitatory pyramidal cells, whereas SST cells predominantly target distal apical dendrites and, in several cortical areas, exert their influence over a broader spatial footprint~\cite{adesnik2012, urbanciecko2016, Kato2017, Aponte2021, kullander2021cortical}. Cortical inhibition is also frequently directed at other interneurons, forming stereotyped disinhibitory circuit motifs: PV cells strongly inhibit one another but provide comparatively little inhibition to other interneuron classes, whereas SST cells largely avoid inhibiting each other yet inhibit PV cells and other populations~\cite{pfeffer2013, jiang2015, kullander2021cortical, campagnola2022local}. Crucially, this diversity is functionally relevant for working memory. Putative interneurons in primate prefrontal cortex carry spatially tuned delay activity~\cite{rao1999, constantinidis2002}, and cell-type-specific manipulations in rodent prefrontal cortex show that both PV and SST interneurons shape delay-period activity and memory-guided behavior, in some task phases with subtype-specific effects~\cite{kamigaki2017, abbas2018}. Inhibitory-to-inhibitory signaling has also been proposed to underlie the stable temporal dynamics required for memory maintenance~\cite{kim2021}, and spiking network models with several interneuron classes have suggested a ``division of labor'' among subtypes in sustaining and sharpening tuned persistent activity~\cite{wang2004}. Theoretical work on circuits with multiple interneuron types, however, has concentrated on stabilization and gain modulation in networks without spatial structure~\cite{litwinkumar2016, garciadelmolino2017, bos2025,kumar2023,negron2024,veit2023}, leaving open how PV and SST architecture shapes the existence, stability, and stochastic precision of spatially localized memory states.

In this work, we introduce a three-population stochastic neural field model with distinct E, PV, and SST populations to investigate how inhibitory subtype structure shapes the existence, stability, and stochastic wandering of persistent activity bumps. The model extends the bipartite E/I framework of~\cite{cihak2022distinct} by partitioning inhibition into a relatively local PV population and a broader SST population, with inter-inhibitory connectivity following the canonical motif of PV self-inhibition and SST-to-PV inhibition, with no SST self-inhibition and weak or absent PV-to-SST coupling, as reported in cortical circuit measurements~\cite{pfeffer2013, campagnola2022local}. In Section~\ref{sec:model}, we formulate the stochastic model and specify the population-specific connectivity architecture. In Section~\ref{sec:stationary}, we construct stationary bump solutions of the deterministic system, and in Section~\ref{sec:linear_stability}, we derive their linear stability by reducing perturbations of the bump interfaces into distinct shifting and scaling modes. Section~\ref{sec:deterministic_effects} uses this stability framework to determine how population thresholds, inhibitory timescales, connection strengths, and spatial scales reshape the stability landscape and the mechanisms by which bumps lose stability. In Section~\ref{sec:wandering}, we return to the stochastic system and derive an effective diffusion coefficient for noise-driven wandering of stable bumps. We conclude with a discussion of the implications and limitations of the model and directions for extending the framework to richer inhibitory circuits.
\section{Model Framework}
\label{sec:model}
To investigate how distinct E, PV, and SST neural populations could influence bump attractor dynamics we utilize a stochastic neural field framework. Extending the bipartite E/I model detailed by \cite{cihak2022distinct}, we construct a system of coupled, nonlinear integro-differential equations to capture the spatial interactions among these three distinct subpopulations. We organize the feature space by the neural populations' angular positional preference, $x\in[-180,180]$, which we take to be in units of degrees. Neural activity feedback from population $n$ to population $m$ at time $t$ are modeled through a convolution of the synaptic weight kernel with the presynaptic firing rates, defined as $(a \ast b)(x) = \int_{-180^\circ}^{180^\circ} a(x - y)b(y)dy$. Incorporating additive stochastic noise, representative of small fluctuations in neuron firing, the governing equations for our neural field are as follows:
\begin{subequations} \label{eq: fullmodel}
    \begin{align}\label{eq: model}
    \tau_edu_{e}(x,t) &= \Big[ -u_{e}(x,t) + w_{ee} \ast f_e(u_{e}) - w_{ep} \ast f_p(u_{p}) - w_{es} \ast f_s(u_{s}) \Big] dt + \epsilon^{1/2} dW_e, \\[1em]
    \tau_p du_{p}(x,t) &= \Big[ -u_{p}(x,t) + w_{pe} \ast f_e(u_{e}) - w_{pp} \ast f_p(u_{p}) - w_{ps} \ast f_s(u_{s}) \Big] dt + \epsilon^{1/2} dW_p, \\[1em]
    \tau_s du_{s}(x,t) &= \Big[ -u_{s}(x,t) + w_{se} \ast f_e(u_{e}) - w_{sp} \ast f_p(u_{p}) - w_{ss} \ast f_s(u_{s}) \Big] dt + \epsilon^{1/2} dW_s.
\end{align}
\end{subequations}

The synaptic input profiles for the E, PV, and SST populations at a specific spatial location $x$ and time $t$ are represented by $u_e(x,t)$, $u_p(x,t)$, and $u_s(x,t)$, respectively. Taking the excitatory population timescale as the reference, we set $\tau_e=1$, corresponding to 10 ms in dimensional time.

The translation of synaptic input $u$ to a firing rate is governed by a sigmoidal function, $f_n(u_n) = \frac{1}{1+e^{-\eta(u_n-\theta_n)}}$, characterized by a gain $\eta$ and a firing threshold $\theta_n$. To facilitate mathematical tractability, we take the high-gain limit ($\eta \to \infty$), which reduces the firing rate to a Heaviside step function:
\begin{equation}\label{eq:heaviside}
    f_n(u_n) = H(u_n - \theta_n) = 
    \begin{cases} 
        1, & u_n - \theta_n \geq 0 \\
        0, & u_n - \theta_n < 0
    \end{cases}.
\end{equation}
Each subpopulation is assigned a potentially distinct firing threshold ($\theta_e$, $\theta_p$, and $\theta_s$). For the Heaviside nonlinearity, localized bump states can be characterized exactly through the interfaces at which $u_n(x,t)=\theta_n$, for $n\in\{e,p,s\}$~\cite{cihak2022distinct,kilpatrick2013wandering}.

\begin{figure}
    \centering
    \includegraphics[width=0.95\linewidth]{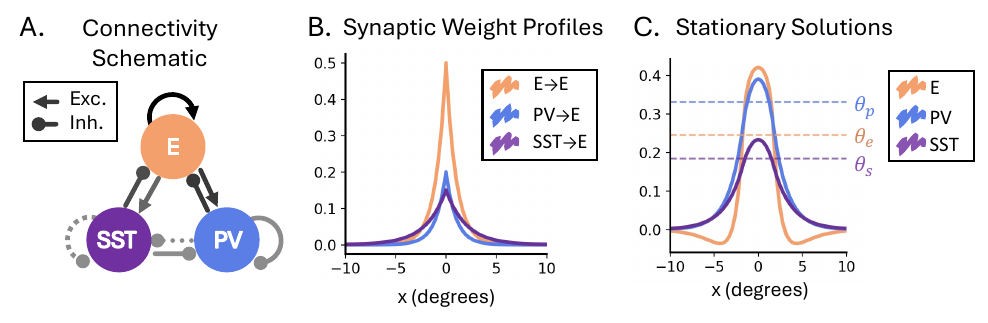}
    \caption{A. Schematic of synaptic connectivity between the E, SST, and PV neuron populations. Dotted lines denote connections set to zero in the baseline parameter set and subsequently varied in Section~\ref{sec:weakIConnect}. B. The spatial parameter $x$ represents the peak orientation stimulus preference of a neuron subpopulation. The synaptic weights were taken to exponentially decaying as neuron stimulus preferences differ, yielding a wizard hat profile. C. Sample stationary bump profile solutions (see section 3 for calculations) for halfwidths $a_e=1.75,a_p=1.25,a_s=1.5$, corresponding to firing thresholds $\theta_e\approx0.245436, \theta_p\approx0.330337, \theta_s\approx0.184118$. For a visualization of the halfwidths see Figure 2.}
    \label{fig:schematics}
\end{figure}

Synaptic coupling strength decays with distance in feature space between the pre- and post-synaptic coordinates, $x$ and $y$. We model these distance-dependent interactions using symmetric exponential weight profiles $w_{mn}(x - y)$, read as an interaction from population $n$ to $m$, for any combination of populations $m,n \in \{e, p, s\}$ on the ring:
\begin{equation}
    w_{mn}(x-y) = A_{mn} e^{-\frac{\min\{|x-y|,360-|x-y|\}}{\sigma_{mn}}}.
\end{equation}
Sample weight kernels, and the corresponding stationary bump profiles derived in Section~\ref{sec:stationary}, are shown in Fig.~\ref{fig:schematics}.

In this formulation, $A_{mn}$ and $\sigma_{mn}$ are nonnegative parameters that determine the peak synaptic strength and spatial decay scale, respectively. We use a baseline parameter set chosen to capture qualitative differences among the three populations rather than to fit a particular cortical circuit. We fix recurrent excitatory coupling at $A_{ee}=0.6$ and $\sigma_{ee}=1$, with outgoing inhibitory connections taken to be weaker than recurrent excitation, reflecting the commonly observed rough ratio of $80\%$ excitatory to $20\%$ inhibitory neurons~\cite{abeles1991corticonics}. To distinguish the spatial roles of the inhibitory subtypes, PV-mediated inhibition is taken to be relatively local, whereas SST-mediated inhibition has a broader spatial footprint~\cite{Kato2017,Aponte2021,campagnola2022local}. Unless otherwise stated, the baseline parameter values used throughout are given in Table~\ref{tab:model_params}.

As a baseline connectivity architecture, we neglect several connections reported to have low connection probability in cortical circuits \cite{campagnola2022local}, setting $A_{sp}=A_{ss}=0$ to zero. These zero values should be understood as modeling approximations rather than an strict biological observation that the corresponding connections are entirely absent. In Section~\ref{sec:weakIConnect}, we relax this approximation to determine how weak inhibitory connections alter bump stability.

To model fluctuations in synaptic input, we incorporate spatially correlated, temporally white additive noise $dW_n(x,t)$, for $n\in\{e,p,s\}$, scaled by a small parameter $0<\epsilon\ll 1$. Each $W_n$ is a mean-zero Wiener process with covariance $\left\langle dW_n(x,t)dW_n(y,s)\right\rangle= C_{nn}(x-y)\delta(t-s)\,dt\,ds$, where where $\delta(z)$ is the Dirac delta distribution. We generate the spatial covariance by filtering spatially white noise $d\Upsilon_n(x,t)$ with the von Mises kernel, $F(x)=\frac{e^{\kappa\cos(x\pi/180)}}{2\pi I_0(\kappa)}$, where $\kappa$ controls the spatial correlation scale and $I_0$ is the modified Bessel function of the first kind of order zero. Specifically, $dW_n(x,t)= \int_{-180^\circ}^{180^\circ}F(x-z)\,d\Upsilon_n(z,t)\,dz$,
so that
\begin{equation*}
    C_{nn}(x-y) = \int_{-180^\circ}^{180^\circ} F(x-z)F(y-z)\,dz,
\end{equation*}
as shown previously in \cite{kilpatrick2013wandering}. Thus, fluctuations are correlated across nearby locations within each
population. Unless otherwise stated, we assume that the noise sources driving
the E, PV, and SST populations are independent, so that $\left\langle dW_m(x,t)dW_n(y,s)\right\rangle=0$, for $m\neq n$. This assumption allows the contributions of fluctuations in each population to the effective diffusion of the coupled bump to be identified separately; correlated noise across populations can be incorporated by retaining the corresponding cross-covariance terms (see ~\cite{cihak2022distinct} for an example of incorporating correlated noise).

\begin{table}[ht]
    \footnotesize
    \caption{Default parameter values used in the computational model.}
    \label{tab:model_params}
    \centering
    \renewcommand{\arraystretch}{1.1} 
    \resizebox{\textwidth}{!}{
    \begin{tabular}{llllllll}
        \hline\noalign{\smallskip}
        Parameter  & $A_{ee}$ & $A_{ep}$& $A_{pe}$ & $A_{es}$ & $A_{se}$ & $A_{pp}, A_{ps}$ & $A_{sp}, A_{ss}$ \\
        Definition & E$\to$E strength & P$\to$E strength & E$\to$P strength & S$\to$E strength & E$\to$S strength & P$\to$P, S$\to$P strength & P$\to$S, S$\to$S strength \\
        Value      & 0.6 & 0.15 & 0.2 & 0.15 & 0.1 & 0.01 & 0.0 \\[2ex]
        Parameter  & \multicolumn{2}{l}{$\sigma_{ee},  \sigma_{pp}$} & $\sigma_{ep}$ & $\sigma_{pe}$ & \multicolumn{3}{l}{$\sigma_{es}, \sigma_{ps}, \sigma_{se}, \sigma_{sp}, \sigma_{ss}$} \\
        Definition & \multicolumn{2}{l}{E$\to$E,  P$\to$P spatial scales} & P$\to$E scale & E$\to$P scale & \multicolumn{3}{l}{Other spatial scales} \\
        Value      & \multicolumn{2}{l}{1.0} & 1.5 & 1.6& \multicolumn{3}{l}{2.0} \\[2ex]
        Parameter  & \multicolumn{3}{l}{$\tau_e,\tau_p, \tau_s$} & $\epsilon$ & $\kappa$ & \multicolumn{2}{l}{}\\
        Definition & \multicolumn{3}{l}{E,PV, SST time constants} & Noise amplitude & Noise spatial scale & \multicolumn{2}{l}{} \\
        Value      & \multicolumn{3}{l}{1.0} & 0.0001 & 1017.0 & \multicolumn{2}{l}{} \\
        \noalign{\smallskip}\hline
    \end{tabular}
    }
\end{table}

Previous studies of neural fields have established conditions under which recurrent excitation balanced by inhibition supports stationary, localized activity ``bumps''~\cite{cihak2022distinct, kilpatrick2013wandering, goldman1995cellular,blomquist2005localized}. In the present model, recurrent E-E coupling provides local excitation, while inhibition from the PV and SST populations acts over distinct spatial scales. This three-population framework allows us to investigate how the balance between relatively local and broad inhibition shapes bump formation and stability. We begin by setting $\epsilon=0$ and characterizing stationary bump solutions of the resulting deterministic system and their linear stability. This analysis identifies regimes in which persistent bump states can be maintained and provides the foundation for the stochastic analysis in Section~\ref{sec:wandering}, where we return to $\epsilon>0$ and quantify the noise-driven wandering of stable bumps.

\section{Stationary Solutions}
\label{sec:stationary}

We begin by characterizing stationary, localized activity profiles, or bump solutions, of the deterministic system. Setting $\epsilon=0$, we seek time-independent solutions $u_n(x,t)=U_n(x)$ for $n\in\{e,p,s\}$. With the Heaviside firing-rate functions defined in Eq.~\ref{eq:heaviside}, the stationary system becomes
\begin{align}\label{eq:stationary}
    U_n(x) &= w_{ne}(x)*H(U_{e}(x)-\theta_e)- w_{np}(x)*H(U_{p}(x)-\theta_p) - w_{ns}(x)*H(U_{s}(x)-\theta_s),
\end{align}
with $n\in\{e,p,s\}$.

\begin{figure}
    \centering
    \includegraphics[width=0.95\linewidth]{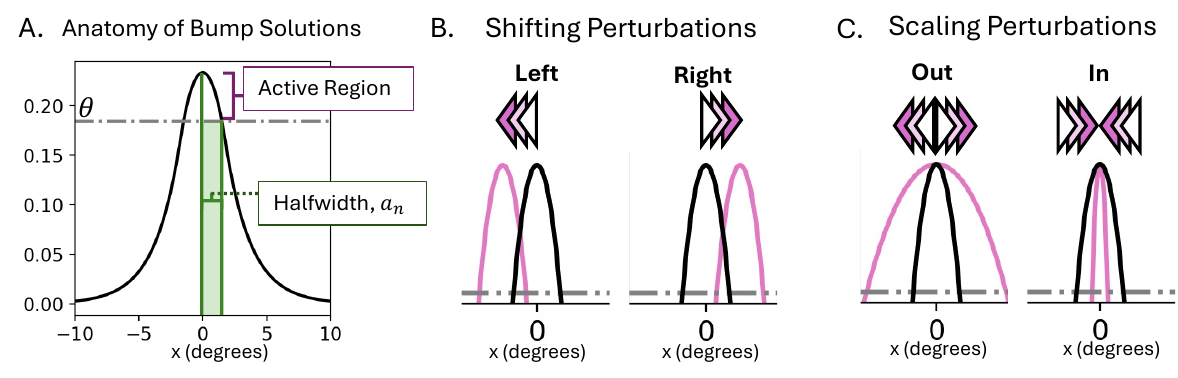}
    \caption{A. Schematic of a bump profile's superthreshold active region and halfwidth $a$ (the distance from the center to the threshold crossing). B. Visualization of ``shift'' perturbations which act to shift the bump from its initial position (black) to the left or right resulting in the final profile (pink). The arrows above denote the shift of the threshold crsossings. C. Visualization of ``scale'' perturbations which act to shift threshold crossings of the bump from its initial position (black) outward or inward resulting in the final profile (pink). The arrows above denote the outward/inward shift of threshold crossings.}
    \label{fig:anatomyPerturbations}
\end{figure}

We restrict attention to stationary solutions for which each population has a single contiguous active region. Specifically, we assume 
\begin{equation*}
    U_n(x) > \theta_n,\quad x \in (-a_n,a_n)
\end{equation*}
with threshold crossings at $x = \pm a_n$, for $n \in \{e,p,s\}$. Here, $a_n$ denotes the half-width of the active region for population $n$ (Fig.~\ref{fig:anatomyPerturbations}A).

As in prior studies \cite{cihak2022distinct,kilpatrick2013wandering,amari1977dynamics}, we expect the system to exhibit translational invariance. Thus, for simplicity and without loss of generality, we pin the bump center to the origin ($x=0$). Restricting attention to symmetric bump solutions, the even synaptic kernels imply $U_n(-x) = U_n(x)$ for $n\in\{e,p,s\}$. With the Heaviside nonlinearity, the assumed active regions $(-a_n,a_n)$, and the symmetry of the stationary profiles, Eq.~\ref{eq:stationary} reduces to
\begin{align}
    U_n(x) &= \int_{-a_e}^{a_e}w_{ne}(x-y)dy - \int_{-a_p}^{a_p}w_{np}(x-y)dy - \int_{-a_s}^{a_s}w_{ns}(x-y)dy,
\end{align}
where $n\in\{e,p,s\}$.
For the symmetric exponential kernel $w_{mn}(x)=A_{mn}e^{\frac{-|x|}{\sigma_{mn}}}$, the convolution over an active interval with half-width $c \in \{a_e, a_p, a_s\}$ can be evaluated explicitly as
\begin{align}
    \int_{-c}^{c}w_{mn}(x-y)dy = \begin{cases}
        2A_{mn}\sigma_{mn}e^{\frac{-x}{\sigma_{mn}}}\sinh\left(\frac{c}{\sigma_{mn}}\right) & \text{if } x>c, \\
        2A_{mn}\sigma_{mn}\left(1-e^{\frac{-c}{\sigma_{mn}}}\cosh\left(\frac{x}{\sigma_{mn}}\right)\right) & \text{if } |x|<c, \\ 
        2A_{mn}\sigma_{mn}e^{\frac{x}{\sigma_{mn}}}\sinh\left(\frac{c}{\sigma_{mn}}\right) & \text{if } x<-c.
    \end{cases}
\end{align}

A necessary self-consistency requirement for these bump solutions is that the synaptic input matches the firing threshold at the interfaces of the active regions; that is, $\theta_n = U_n(\pm a_n)$ for $n\in\{e,p,s\}$ or
\begin{align*}
    \theta_n &= U_n(a_n)=\int_{-a_e}^{a_e}w_{ne}(a_n-y)dy - \int_{-a_p}^{a_p}w_{np}(a_n-y)dy - \int_{-a_s}^{a_s}w_{ns}(a_n-y)dy,
\end{align*}
for $n\in\{e,p,s\}$. Because the explicit form of each integral depends on whether $a_n$ lies inside or outside the other active regions, the resulting threshold conditions are piecewise with respect to the ordering of $a_e$, $a_p$, and $a_s$:
\begin{subequations}\label{eqs: piecewise_stationary}
\begin{align}
    \theta_e = 2A_{ee}\sigma_{ee}\mathcal{S}_{ee} - \begin{cases} 
        2A_{ep}\sigma_{ep}\mathcal{S}_{ep} + 2A_{es}\sigma_{es}\mathcal{S}_{es} & \text{if } a_e \geq \max\{a_p,a_s\} \\
        2A_{ep}\sigma_{ep}\mathcal{C}_{ep} + 2A_{es}\sigma_{es}\mathcal{S}_{es} & \text{if } a_s \leq a_e < a_p \\
        2A_{ep}\sigma_{ep}\mathcal{S}_{ep} + 2A_{es}\sigma_{es}\mathcal{C}_{es} & \text{if } a_p \leq a_e < a_s \\
        2A_{ep}\sigma_{ep}\mathcal{C}_{ep} + 2A_{es}\sigma_{es}\mathcal{C}_{es} & \text{if } a_e < \min\{a_p,a_s\} 
       \end{cases}
\end{align}
\begin{align}
    \theta_p = -2A_{pp}\sigma_{pp}\mathcal{S}_{pp} + \begin{cases} 
        2A_{pe}\sigma_{pe}\mathcal{S}_{pe} - 2A_{ps}\sigma_{ps}\mathcal{S}_{ps} & \text{if } a_p \ge \max\{a_e,a_s\} \\ 
        2A_{pe}\sigma_{pe}\mathcal{C}_{pe} - 2A_{ps}\sigma_{ps}\mathcal{S}_{ps} & \text{if } a_s \leq a_p < a_e \\
        2A_{pe}\sigma_{pe}\mathcal{S}_{pe} - 2A_{ps}\sigma_{ps}\mathcal{C}_{ps} & \text{if } a_e \le a_p < a_s \\
        2A_{pe}\sigma_{pe}\mathcal{C}_{pe} - 2A_{ps}\sigma_{ps}\mathcal{C}_{ps} & \text{if } a_p < \min\{a_e,a_s\}
       \end{cases}
\end{align}
\begin{align}
    \theta_s = -2A_{ss}\sigma_{ss}\mathcal{S}_{ss} + \begin{cases} 
        2A_{se}\sigma_{se}\mathcal{S}_{se} - 2A_{sp}\sigma_{sp}\mathcal{S}_{sp} & \text{if } a_s \ge \max\{a_e, a_p\} \\
        2A_{se}\sigma_{se}\mathcal{C}_{se} - 2A_{sp}\sigma_{sp}\mathcal{S}_{sp} & \text{if } a_p \leq a_s < a_e \\
        2A_{se}\sigma_{se}\mathcal{S}_{se} - 2A_{sp}\sigma_{sp}\mathcal{C}_{sp} & \text{if } a_e \le a_s < a_p \\
        2A_{se}\sigma_{se}\mathcal{C}_{se} - 2A_{sp}\sigma_{sp}\mathcal{C}_{sp} & \text{if } a_s < \min\{a_e,a_p\} 
       \end{cases}
\end{align}
\end{subequations}
where we define the functions
\begin{align*}
    \mathcal{S}_{mn} &= e^{-\frac{a_m}{\sigma_{mn}}}\sinh\left(\frac{a_n}{\sigma_{mn}}\right), \quad
    \mathcal{C}_{mn} = 1-e^{-\frac{a_n}{\sigma_{mn}}}\cosh\left(\frac{a_m}{\sigma_{mn}}\right),
\end{align*}
which correspond respectively to evaluation of the convolution at an interface outside or inside the presynaptic active region.

The stationary solutions can be obtained numerically in two complementary ways. We may prescribe the half-widths $(a_e, a_p, a_s)$, compute the corresponding thresholds from Eq.~\ref{eqs: piecewise_stationary}, and then verify that the resulting profiles satisfy the assumed single-active-region structure. Alternatively, for prescribed thresholds $(\theta_e, \theta_p, \theta_s)$, we solve the nonlinear system for the half-widths using a Levenberg-Marquardt algorithm and similarly verify that resulting profiles satisfy the assumed single-active-region structure. Both methods were used in different figures and which parameters were fixed are noted in the captions.

With the stationary solutions determined, we next examine their linear stability to small perturbations localized at the interfaces of the active regions. As illustrated in Fig.~\ref{fig:anatomyPerturbations}B,C, these perturbations can be separated into two symmetry classes: shifting perturbations, which translate the profiles to the left or right, and scaling perturbations, which expand or contract the active regions. This distinction will allow the linear stability problem to be decomposed into separate shift and scale modes in Section~\ref{sec:linear_stability}.

\section{Linear Stability Analysis}
\label{sec:linear_stability}
To characterize the linear stability of the stationary bump solutions, we consider small spatiotemporal perturbations about the equilibrium profiles. For $n\in\{e,p,s\}$, let
\begin{align*}
    u_n(x,t) &= U_n(x) + \nu \Psi_n(x,t),
\end{align*}
where $0 < \nu \ll 1$ captures the perturbation amplitude and $\Psi_n$ denotes the corresponding perturbation function. We linearize the firing-rate function about the stationary solution using the first-order expansion
\begin{equation*}
    f(\alpha(x) + \nu \beta(x,t)) \approx f(\alpha(x)) + f'(\alpha(x))(\nu \beta(x,t)).
\end{equation*}
Substituting the perturbed states into the deterministic system and retaining terms of $\mathcal{O}(\nu)$ gives the linearized equations
\begin{align}\label{eq:evol_perturbations}
    \tau_n\frac{\partial \Psi_n}{\partial t} &= -\Psi_n(x,t) + w_{ne} \ast [f'(U_e)\Psi_e] - w_{np} \ast [f'(U_p) \Psi_p] - w_{ns} \ast [f'(U_s)\Psi_s]
\end{align}
for $n\in\{e,p,s\}$. We assume the perturbations are separable in space and time and take the form $\Psi_n(x,t) = \psi_n(x)e^{\lambda t}$. Substituting and factoring out the temporal components in  Eq.~\ref{eq:evol_perturbations} reduces the system to
\begin{align}\label{eq:separated_perturb}
    (\tau_n\lambda + 1)\psi_n(x) &= w_{ne} \ast [f'(U_e)\psi_e] - w_{np} \ast [f'(U_p)\psi_p] - w_{ns} \ast [f'(U_s)\psi_s].
\end{align}
For the Heaviside firing-rate function, $f(u) = H(u-\theta)$, the derivative is interpreted distributionally. Using 
\begin{equation*}
    f'(U_n(x)) = \delta(U_n(x) - \theta_n), 
\end{equation*}
and the threshold crossings $U_n(\pm a_n) = \theta_n$, we obtain
\begin{equation}\label{eq:fprime}
    f'(U_n(x)) = \frac{1}{|U_n'(a_n)|} \Big( \delta(x-a_n) + \delta(x+a_n) \Big).
\end{equation}
Applying Eq.~\ref{eq:fprime} to Eq.~\ref{eq:separated_perturb} localizes the convolution terms to the threshold interfaces. To simplify the resulting expression, we define
\begin{equation*}
    \psi_n^\pm := \psi_n(\pm a_n), \quad n\in\{e,p,s\}.
\end{equation*}
The eigenvalue problem then becomes
\begin{align} \label{eq:sixEquations}
    (\tau_n\lambda + 1)\psi_n(x) &= \frac{A_{ne}}{|U_e'(a_e)|}\left[\exp\left(-\frac{|x-a_e|}{\sigma_{ne}}\right)\psi_e^+ + \exp\left(-\frac{|x+a_e|}{\sigma_{ne}}\right)\psi_e^-\right]  \\
    &\quad - \frac{A_{np}}{|U_p'(a_p)|}\left[\exp\left(-\frac{|x-a_p|}{\sigma_{np}}\right)\psi_p^+ + \exp\left(-\frac{|x+a_p|}{\sigma_{np}}\right)\psi_p^-\right] \notag \\
    &\quad - \frac{A_{ns}}{|U_s'(a_s)|}\left[\exp\left(-\frac{|x-a_s|}{\sigma_{ns}}\right)\psi_s^+ + \exp\left(-\frac{|x+a_s|}{\sigma_{ns}}\right)\psi_s^-\right],\notag
\end{align}
for $n\in\{e,p,s\}$. We next evaluate Eq.~\ref{eq:sixEquations} at the six threshold interfaces $x = \pm a_n$, $n\in\{e,p,s\}$. For compactness, we define
\begin{align*}
    \chi_e = 1, \chi_p = \chi_s = -1,
\end{align*}
and the interaction coefficients as
\begin{align*}
    K_{mn}^\pm = \chi_n \frac{A_{mn}}{|U'_n(a_n)|}\exp\left(-\frac{|a_m\pm a_n|}{\sigma_{mn}}\right), \quad m,n\in\{e,p,s\}.
\end{align*}
Here, $K^-_{mn}$ describes the interaction between interfaces on the same side of the bump, while $K^+_{mn}$ corresponds to coupling between opposite-side interfaces. This evaluation yields six equations of the form
\begin{align*}
    0 = -(1+\tau_m \lambda)\psi_m^{\pm} + \sum_{n\in\{e,p,s\}} \left(K^-_{mn}\psi_n^\pm+K^+_{mn}\psi_n^\mp\right), \quad m \in \{e,p,s\}.
\end{align*}
These six equations define a closed six-dimensional eigenvalue problem for the interface perturbations $\psi_n^\pm$. Although the corresponding characteristic equation is sixth order in $\lambda$, reflection symmetry of the stationary bump solutions allows the perturbations to be decomposed into odd and even modes, reducing the problem to two independent three-dimensional eigenvalue problems. These correspond, respectively, to the shifting and scaling perturbations illustrated in Fig.~\ref{fig:anatomyPerturbations}B,C.

\subsection{Shifting Perturbations (Translation)}\label{sec:shift}
Shifting perturbations are odd with respect to the common center of symmetry, so that $\psi_n^+ = -\psi_n^-$ for $n\in\{e,p,s\}$, as illustrated in Fig.~\ref{fig:anatomyPerturbations}B. The odd subspace contains the translational mode associated with shifting the stationary bump solution along the feature space. Substituting this symmetry into the six-dimensional interface system reduces the eigenvalue problem to
\begin{align*}
    {\bf 0} =\mathbf{M}_{\text{shift}}(\lambda)\begin{pmatrix}
        \psi_e^+ \\ \psi_p^+ \\ \psi_s^+
    \end{pmatrix},
\end{align*}
with
\begin{equation} \label{eq:shiftM}
\mathbf{M}_{\text{shift}}(\lambda) = 
    \begin{pmatrix} 
        K_{ee}^--K_{ee}^+ - 1- \lambda & K_{ep}^--K_{ep}^+ & K_{es}^--K_{es}^+ \\
        K_{pe}^--K_{pe}^+ & K_{pp}^--K_{pp}^+ - 1- \tau_p\lambda & K_{ps}^--K_{ps}^+ \\
        K_{se}^--K_{se}^+ & K_{sp}^--K_{sp}^+ & K_{ss}^--K_{ss}^+- 1- \tau_s\lambda \\
    \end{pmatrix}.
\end{equation}
Continuous translational invariance guarantees a zero eigenvalue, $\lambda=0$, in the shifting subspace. The associated eigenmode is neutrally stable and corresponds to translation of the stationary bump along the feature space. Stability with respect to shifting perturbations requires all remaining eigenvalues of $\mathbf{M}_{\text{shift}}(\lambda)$ to have negative real part.

\subsection{Scaling Perturbations (Expansion/Contraction)}
Scaling perturbations are even with respect to the common center of symmetry, meaning $\psi_n^+ = \psi_n^-$ for $n\in\{e,p,s\}$. These perturbations correspond to symmetric expansion or contraction of the active regions, as illustrated in Fig.~\ref{fig:anatomyPerturbations}C. Substituting this symmetry into the six-dimensional interface system reduces the eigenvalue problem to
\begin{align*}
    {\bf 0} = \mathbf{M}_{\text{scale}}(\lambda)\begin{pmatrix}
        \psi_e^+ \\ \psi_p^+ \\ \psi_s^+
    \end{pmatrix},
\end{align*}
where
\begin{equation}
    \label{eq:scaleM}
    \mathbf{M}_{\text{scale}}(\lambda) = 
    \begin{pmatrix} 
        K_{ee}^-+K_{ee}^+ - 1- \lambda & K_{ep}^-+K_{ep}^+ & K_{es}^-+K_{es}^+ \\
        K_{pe}^-+K_{pe}^+ & K_{pp}^-+K_{pp}^+ - 1- \tau_p\lambda & K_{ps}^-+K_{ps}^+ \\
        K_{se}^-+K_{se}^+ & K_{sp}^-+K_{sp}^+ & K_{ss}^-+K_{ss}^+- 1- \tau_s\lambda \\
    \end{pmatrix}.
\end{equation}
Stability with respect to scaling perturbations requires all eigenvalues of $\mathbf{M}_{\text{scale}}(\lambda)$ to have negative real part.

\subsection{Stability Classification}\label{sec:stab_class}
The preceding analysis separates perturbations into shifting and scaling modes, with stability requiring decay of all non-translational modes in both subspaces. To compare the relative stability of these modes as model parameters are varied, let
$\{\lambda_j^{\mathrm{shift}}\}$ and $\{\lambda_j^{\mathrm{scale}}\}$ denote
the eigenvalues associated with the shifting and scaling subspaces, respectively. We define
\begin{equation*}
\Lambda_{\mathrm{shift}} = \max_{\lambda_j^{\mathrm{shift}}\neq 0} \operatorname{Re}\!\left(\lambda_j^{\mathrm{shift}}\right), \qquad \Lambda_{\mathrm{scale}}=\max_j\operatorname{Re}\!\left(\lambda_j^{\mathrm{scale}}\right),
\end{equation*}
where the translational zero eigenvalue is excluded from
$\Lambda_{\mathrm{shift}}$. The overall stability metric is then
\begin{equation*}
\Lambda = \max\left\{\Lambda_{\mathrm{shift}}, \Lambda_{\mathrm{scale}} \right\}.
\end{equation*}
A stationary bump is linearly stable when $\Lambda<0$. When $\Lambda>0$,
comparison of $\Lambda_{\mathrm{shift}}$ and $\Lambda_{\mathrm{scale}}$
identifies whether the dominant instability lies in the shifting or scaling subspace. In Section~\ref{sec:deterministic_effects}, we use these quantities to map stability across parameter space and distinguish the mechanisms by which changes in inhibitory subtype properties destabilize the bump.

\section{Deterministic effects of inhibitory subtype structure}
\label{sec:deterministic_effects}
Using the stability metrics defined in Section~\ref{sec:stab_class}, we now examine how inhibitory subtype properties shape the existence and stability of stationary bump solutions. We begin by varying the firing thresholds of the three populations to map the region of stable solutions in the full three-dimensional threshold space 
and then use lower-dimensional parameter slices to isolate the effects of population thresholds and inhibitory timescales.

In Section~\ref{sec:bifurication}, we characterize how differences in PV and SST recruitment and timescale can alter both bump stability and the mechanism by which stability is lost. We then turn to the connectivity architecture itself in Section~\ref{sec:weakIConnect}, varying inhibitory connection strengths and spatial scales to determine how these features reshape the stability landscape. Together, these analyses highlight complementary roles for localized PV and broader SST inhibition, while also showing that even weak connections among inhibitory populations can have substantial, and potentially nonmonotonic, effects on bump stability.

\subsection{Stable Parameter Region Boundaries}
\label{sec:bifurication}
\begin{figure}
    \centering
    \includegraphics[width=\linewidth]{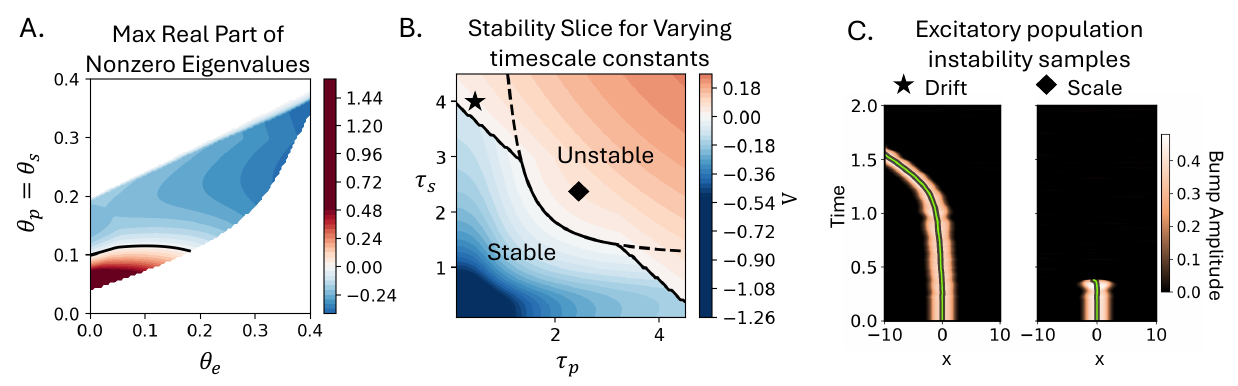}
    \caption{A. A parameter slice where the thresholds of the inhibitory populations are set equal and parameters are as in Table \ref{tab:model_params}. B. Increasing the timescales of PV and SST neurons results in eventual destabilization. Interestingly, we can obtain two types of instability dependent on a comparison of $\Lambda_\text{shift}$ and $\Lambda_\text{scale}$. The solid line marks the separation between stable and unstable space. The dashed line indicates where $\Lambda_\text{scale}$ is positive which will lead to oscillatory or scale instabilities (region marked with a diamond). It is possible for $\Lambda_\text{shift}$ to become positive and cause a drift instability or a traveling wave with $\Lambda_\text{scale}<0$ which would be the regions at the flanks with the upper left marked with a star. Bump solution halfwidths were fixed as $(a_e,a_p,a_s)=(1.75,1.25,1.5)$. C. Samples of each type of instability with the same noisy perturbations. Left: drift instability at $(\tau_p,\tau_s)=(0.5,4)$ can result in traveling waves. Right: A scale or oscillatory instability $(\tau_p,\tau_s)=(2.5,2.5)$ results in silencing of the memory.}
    \label{fig:stabilityA}
\end{figure}

We first examine a two-dimensional slice of the threshold space in Fig.~\ref{fig:stabilityA}A by setting the inhibitory thresholds equal, $\theta_s=\theta_p$, while retaining the distinct PV and SST connectivity parameters. For fixed $\theta_e$, increasing the common inhibitory threshold generally moves the system into the stable region. At the boundary marked in black, a complex-conjugate pair of eigenvalues crosses the imaginary axis, indicating a Hopf bifurcation. This organization of the stability boundary is similar to that found in the two-population model of \cite{cihak2022distinct}, providing a useful point of consistency between the models. The additional inhibitory degree of freedom in the present model, however, permits a richer set of instability dynamics, which we examine next by varying the PV and SST timescales. 

Varying the inhibitory timescales, $\tau_p$ and $\tau_s$, we use $\Lambda$ to determine bump stability and $\Lambda_\text{shift}$ and $\Lambda_\text{scale}$ to identify the dominant instability mode (Fig.~\ref{fig:stabilityA}B). When the inhibitory populations respond sufficiently rapidly relative to the excitatory population, the bump remains stable. As either inhibitory timescale increases, however, the system can lose stability through either the shifting or scaling subspace. A dominant shifting instability produces drift or traveling-wave dynamics (Fig.~\ref{fig:stabilityA}C, left), whereas a dominant scaling instability produces oscillatory changes in bump width that can ultimately annihilate the bump (Fig.~\ref{fig:stabilityA}C, right).
\begin{figure}
    \centering
    \includegraphics[width=\linewidth]{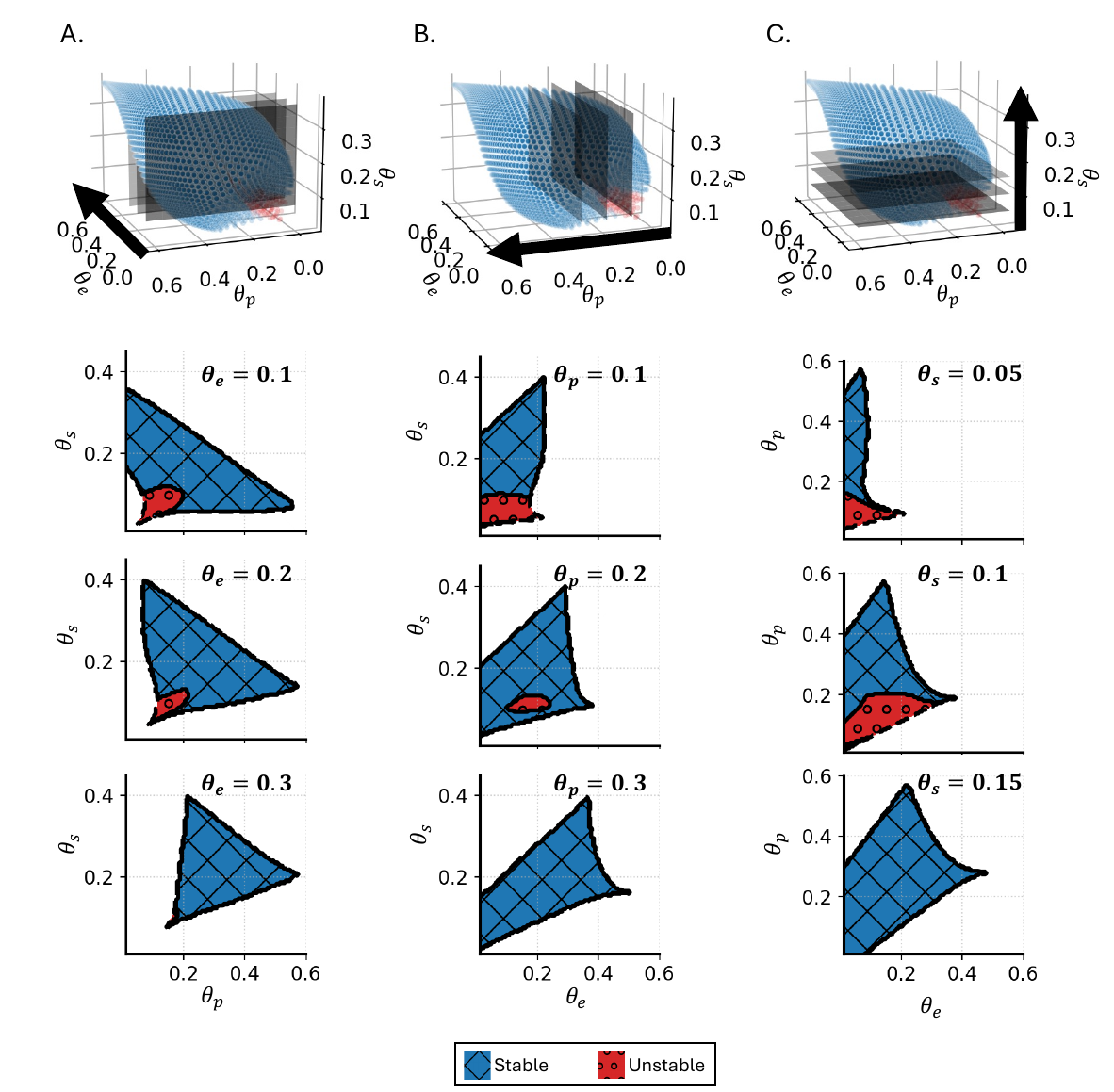}
    \caption{Slices taken across parameter space along different threshold axes. Columns A, B, and C represent traces along axes $\theta_e, \theta_p$ and $\theta_s$ respectively. The arrow in top 3D figure demonstrates the direction of this axis. Below are three panels showing the existence and stability of bumps for increasing thresholds along each axis respectively.}
    \label{fig:stability3d}
\end{figure}

To separate the effects of the population thresholds, we next examine two-dimensional slices through the full threshold space (Fig.~\ref{fig:stability3d}). Increasing the threshold of any population generally reduces the unstable portion of the corresponding slice, but the stability boundaries are not uniformly monotone. For example, at the intermediate value $\theta_p=0.2$, stable solutions occur for relatively low values of $\theta_e$ and $\theta_s$, while an unstable region emerges at intermediate threshold values. This behavior highlights that stability is determined by the balance among the three populations rather than by the threshold of any one population independently.

The threshold slices also reveal a notable difference between the two inhibitory populations. Changes in the SST threshold produce a stronger restriction of the region admitting stationary bumps than comparable changes in the PV threshold (Fig.~\ref{fig:stability3d}B,C). Within our baseline architecture, SST provides the broader component of inhibition, suggesting that recruitment of this longer-range inhibitory population plays an important role in determining the existence of localized bump states. More generally, these results emphasize that both bump existence and stability depend on the balance between the distinct PV and SST contributions to inhibition.
\subsection{Effects of Inhibitory Connectivity Structure}
\label{sec:weakIConnect}
The baseline connectivity architecture used thus far is motivated by cortical connectivity measurements in which several connections among inhibitory populations are relatively weak \cite{campagnola2022local, riedemann2019diversity}. Accordingly, some of these connections were set to zero in the baseline model. Here, we relax this approximation to determine whether weak inhibitory connections, and their corresponding spatial footprint, can have substantial effects on bump stability.

To quantify these effects, we examine bump stability over a fixed set of sampled points
\begin{equation*}
   \mathcal{P}=\left\{q_i=(\tau_{p,i},\tau_{s,i})\right\}_{i=1}^{N}. 
\end{equation*}
in the inhibitory-timescale plane while varying a single connectivity parameter. To track how this stability landscape changes as connectivity parameter is varied, we introduce the following notation. Let $\Lambda(q_i;\alpha)$ be the stability metric defined in Section~\ref{sec:stab_class} at the point $q_i$ for a fixed connectivity value $\alpha$, whenever a valid stationary bump exits. We then define the stable subset
\begin{equation*}
    \mathcal{S}(\alpha) = \left\{ q_i\in\mathcal{P}: \Lambda(q_i;\alpha)<0 \right\},
\end{equation*}
and the corresponding proportion of the sampled plane that is stable as
\begin{equation*}
    P_{\mathrm{stable}}(\alpha) = \frac{|\mathcal{S}(\alpha)|}{N}.
\end{equation*}
We additionally define $\Lambda_{\min}(\alpha)$ as the minimum value of $\Lambda(q_i;\alpha)$ over those sampled points for which it is defined, which measures the most negative value of the stability metric attained within the sampled plane. When no valid stationary bumps (i.e., bumps that are positive and have a single connected active region) are found within the sampled plane, $\Lambda_{\min}(\alpha)$ is left undefined.

\begin{figure}
    \centering
    \includegraphics[width=0.95\linewidth]{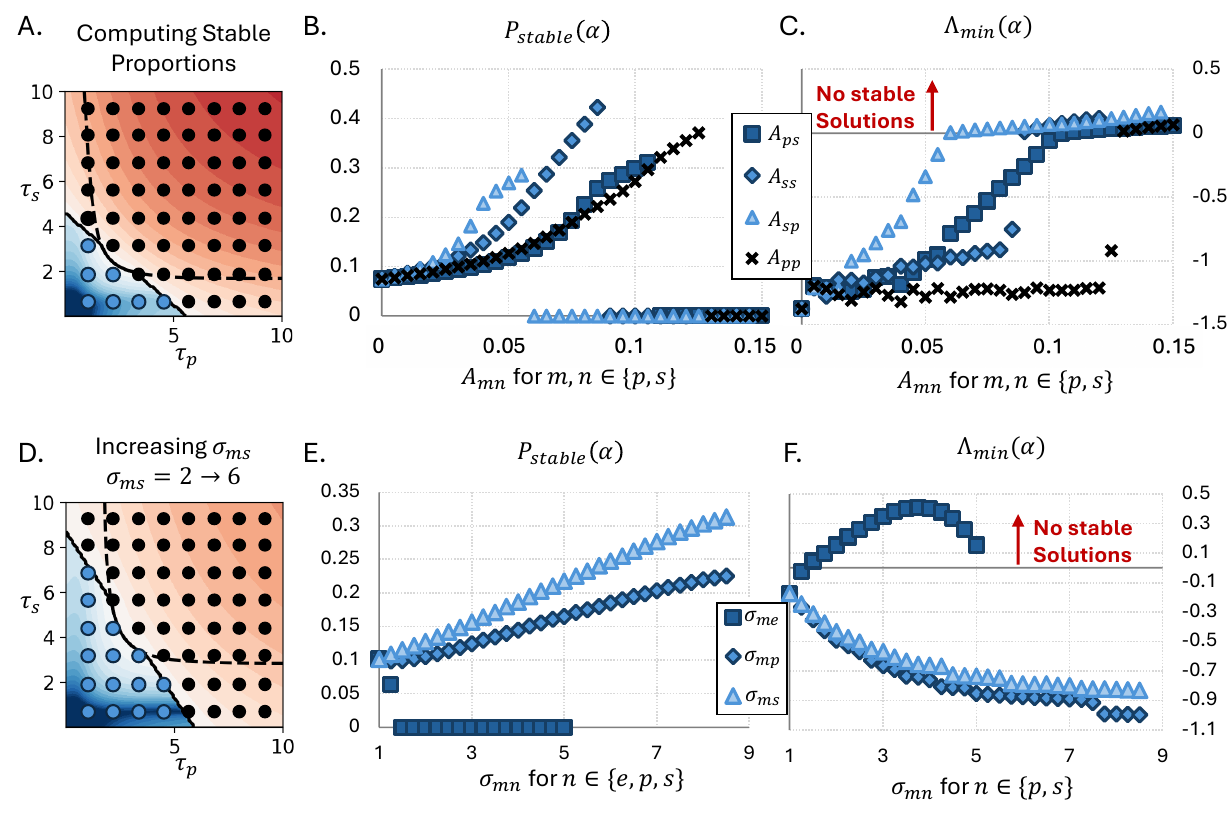}
    \caption{A. The stable proportion of a snapshot of parameter space is determined by summing the number of stable solutions (blue dots) over the total number sampled (10000 samples evenly spaced in $[0.1,10]\times[0.1,10]$). Bumps are classified by the max real part of the solution's nonzero eigenvalues. In this snapshot, $A_{mn}=0.01$ for $m,n\in\{p,s\}$ B. The proportion of the slice that is stable and C. the corresponding minimum max nonzero eigenvalue as individual inhibitory connections are increased (all inhibitory to inhibitory connections are initialized at zero). Values above zero in panel C indicate the whole parameter slice snapshot has become unstable. D,E,F demonstrate similar plots but rather for varying synaptic footprints $\sigma_{mn}$. D. The same parameters as in panel A save for $\sigma_{es}=\sigma_{ps}=\sigma_{ss}$ which were tripled to 6. E \& F. Same as B\& C but with varying synaptic footprints. For $\sigma_{me}$ valid solutions were not found beyond the value of 5 and thus no data points are pictured. For D,E, and F  $A_{mn}=0.01$ for $m,n\in\{p,s\}$.  All parameters are as in Table \ref{tab:model_params} unless otherwise noted or varied in the graph traces. Halfwidths were fixed as $(a_e,a_p,a_s)=(1.75,1.25,1.5)$.}
    \label{fig:proportionStable}
\end{figure}

We first consider changes in inhibitory connection strength $A_{mn}$ for $m,n\in \{p,s\}$ (in other words $\alpha=A_{mn}$). Figure~\ref{fig:proportionStable}A shows a representative stability landscape over the sampled $(\tau_p,\tau_s)$ plane, with each point classified according to whether the corresponding bump solution is stable or unstable. Repeating this over the same sampled plane as the connectivity strength $\alpha$ is varied allows us to track both the proportion of stable points, $P_\text{stable}(\alpha)$, and the minimum stability metric, $\Lambda_\text{min}(\alpha)$. As we increase the strength of individual inhibitory connections while recomputing stability over the same sampled timescale plane, the proportion of stable points initially increases, indicating an expansion of the region supporting stable bumps (Fig.~\ref{fig:proportionStable}B). Beyond a critical connection strength, however, this trend reverses and the stable proportion falls to zero. 

To better understand this nonmonotonic behavior, we examine $\Lambda_{\min}(\alpha)$ over the same connectivity sweeps (Fig.~\ref{fig:proportionStable}C). As the stable region expands, $\Lambda_{\min}(\alpha)$ simultaneously moves toward zero, indicating that even the most strongly attracting point in the sampled plane becomes progressively less stable. With further increases in connectivity strength, $\Lambda_{\min}(\alpha)$ becomes positive, at which point the entire sampled timescale plane is unstable. Thus, weak inhibitory coupling can broaden the region of parameter space supporting stable bumps while simultaneously weakening the strongest linear restoring rate attained within that region. These results are consistent with other modeling studies that found strong inhibitory to inhibitory connections destabilized the bumps~\cite{blomquist2005localized, pinto2001spatially}. Furthermore, they show that even weak inhibitory to inhibitory connections can substantially reshape the stability landscape.

We next apply the same analysis to changes in the spatial scales $\sigma_{mn}$ ($\alpha= \sigma_{mn}$ in this case). To isolate the effect of projection range, we vary together the spatial footprints of projections originating from a given presynaptic population, considering $\sigma_{me}$, $\sigma_{mp}$, and $\sigma_{ms}$ for $m\in\{e,p,s\}$ while holding the remaining parameters fixed. Fig.~\ref{fig:proportionStable}D illustrates a example plane where we broadened the SST projections, which produced a visible expansion of the stable region in the inhibitory-timescale plane. Repeating these spatial-scale sweeps and computing $P_{\mathrm{stable}}$ and $\Lambda_{\min}(\alpha)$ reveals distinct effects of excitatory and inhibitory projection range (Fig.~\ref{fig:proportionStable}E,F). Increasing the spatial footprint of either the PV or SST population increases the proportion of the sampled plane that is stable while simultaneously driving $\Lambda_{\min}(\alpha)$ to more negative values. Thus, in contrast to increasing inhibitory connection strength, broadening inhibitory projections both expands the stable region and strengthens the largest linear restoring rates within it. For excitatory projections, increasing the spatial footprint instead eventually prevents the localized inhibitory populations from supporting the assumed bump structure; beyond this regime, valid single-active-region solutions are no longer obtained. Although the model predicts a stabilizing effect from broadening either inhibitory population, PV neurons are expected to have a more localized spatial reach than SST neurons. These results therefore suggest more generally that the spatial extent of inhibition can strongly influence bump stability, and support a complementary role for localized and longer-range inhibitory connectivity in maintaining persistent activity.

Together, these results demonstrate that features of inhibitory connectivity that may appear secondary based on connection strength alone can substantially alter the existence and linear stability of persistent bump states. Moreover, the changes in $\Lambda$ observed across these parameter sweeps suggest that inhibitory connectivity may influence not only whether a bump is stable, but also how strongly perturbations away from the stationary state decay. We next ask whether these deterministic changes in stability are accompanied by corresponding changes in the stochastic wandering of the memory representation.

\section{Diffusion of Wandering Bumps}
\label{sec:wandering}
\begin{figure}
    \centering
    \includegraphics[width=\linewidth]{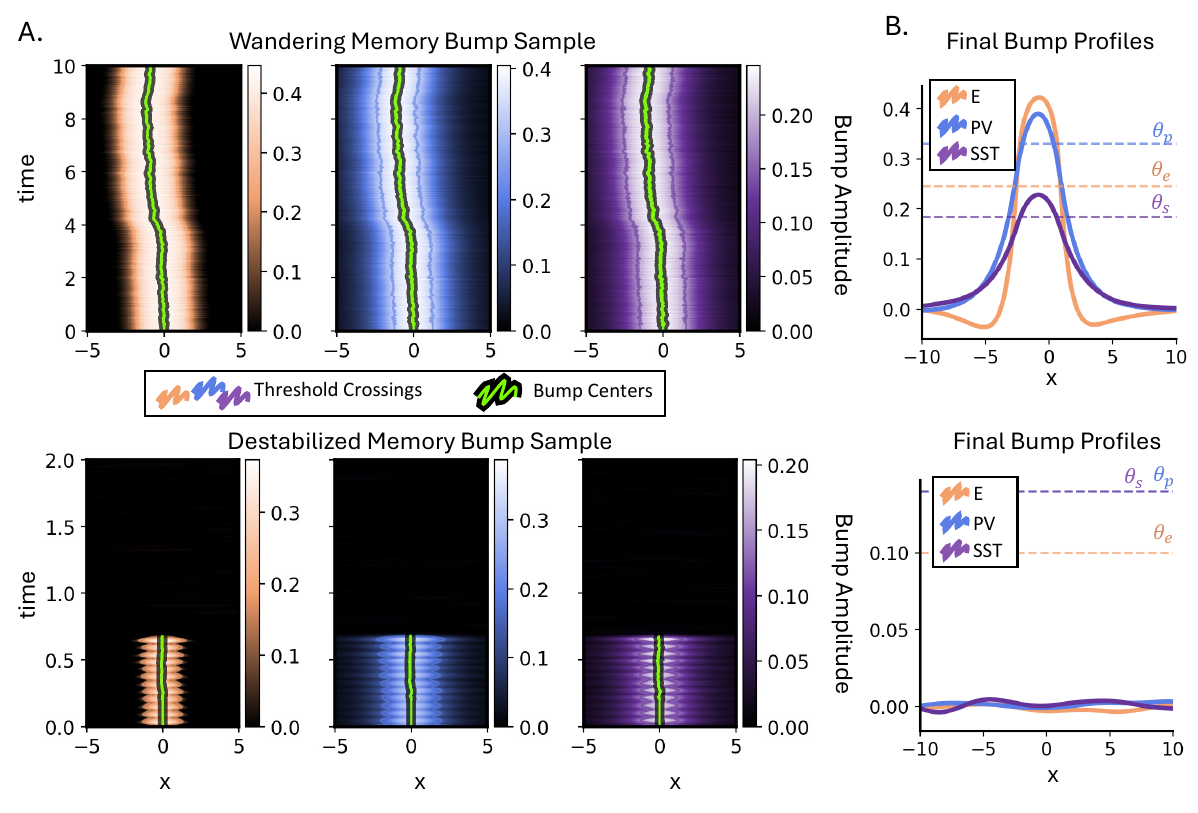}
    \caption{A. Heatmaps of bump solutions in response to noise. Standard parameters with halfwidths $(a_e,a_p,a_s)=(1.75,1.25,1.5)$ is pictured above. The bottom panel shows a solution very close to the boundary of a scaling instability, resulting in the annihilation of the bump. In this case thresholds were set $(\theta_e,\theta_p,\theta_s)=(0.1,0.14,0.14)$. All other parameters are as in Table \ref{tab:model_params}.}
    \label{fig:wanderingsample}
\end{figure}

We first confirm numerically that stable bump solutions exhibit stochastic wandering under weak noise (Fig.~\ref{fig:wanderingsample}). For solutions well within the stable region, the three population profiles remain localized while their common position fluctuates over time (Fig.~\ref{fig:wanderingsample}, upper panels). In contrast, solutions near a stability boundary can be driven toward the deterministic instabilities identified in Section~\ref{sec:deterministic_effects}, resulting in traveling activity or annihilation of the bump (Fig.~\ref{fig:wanderingsample}, lower panels). The variance of the bump position provides a natural measure of this wandering and has been linked to behavioral response errors in delayed-estimation tasks~\cite{compte2000synaptic,wimmer2014bump}. We therefore seek an effective stochastic description of the common bump position $\Delta(t)$ and, in particular, a prediction for its variance $\langle\Delta(t)^2\rangle$ in the small-noise limit.

For a linearly stable bump, all non-translational perturbations decay, while translational invariance leaves the zero eigenvalue identified in Section~\ref{sec:shift}. Thus, under weak stochastic forcing, the leading long-time displacement occurs along the common translational mode of the three population profiles. We denote this displacement by $\Delta(t)$ and assume $d\Delta=\mathcal{O}(\epsilon^{0.5})$, consistent with the noise amplitude in Eq.~\ref{eq: fullmodel}. We then see a weak-noise expansion through first order in $\epsilon^{0.5}$ of the form
\begin{equation*}
    u_n(x,t)= U_n(x-\Delta(t))+\epsilon^{0.5}\Psi_n(x-\Delta(t),t)+\mathcal{O}(\epsilon), \quad n\in \{e,p,s\},
\end{equation*}
where $\Psi_n$ captures small fluctuations transverse to the common translational mode. Substituting this ansatz into Eq.~\ref{eq: fullmodel} and retaining terms through $\mathcal{O}(\epsilon^{0.5})$ yields the following leading-order stochastic perturbation equations
\begin{equation}\label{eq:wandereq1}
    \begin{pmatrix}
        \tau_ed\Psi_e\\
        \tau_pd\Psi_p\\
        \tau_sd\Psi_s
    \end{pmatrix}=\mathcal{L}\begin{pmatrix}
        \Psi_e\\\Psi_p\\\Psi_s
    \end{pmatrix}dt+\epsilon^{-0.5}d\Delta(t)\begin{pmatrix}
        \tau_e U_e'(x)\\\tau_pU_p'(x)\\\tau_sU_s'(x)
    \end{pmatrix}+\begin{pmatrix}
        dW_e\\dW_p\\dW_s
    \end{pmatrix}
\end{equation}
where $\mathcal{L}$ is the linearized deterministic operator about the stationary bump, namely
\begin{equation}
    \mathcal{L}\begin{pmatrix}
        \Psi_e\\\Psi_p\\\Psi_s
    \end{pmatrix}=\begin{pmatrix}
        -\Psi_e+w_{ee}*\left[f'(U_e)\Psi_e\right]-w_{ep}*\left[f'(U_p)\Psi_p\right]-w_{es}*\left[f'(U_s)\Psi_s\right]\\-\Psi_p+w_{pe}*\left[f'(U_e)\Psi_e\right]-w_{pp}*\left[f'(U_p)\Psi_p\right]-w_{ps}*\left[f'(U_s)\Psi_s\right]\\-\Psi_s+w_{se}*\left[f'(U_e)\Psi_e\right]-w_{sp}*\left[f'(U_p)\Psi_p\right]-w_{ss}*\left[f'(U_s)\Psi_s\right]
    \end{pmatrix}.
\end{equation}
Similar to \cite{cihak2022distinct}, we find that due to the translation symmetry of the system the nullspace of the operator is spanned by the vector $(U_e',U_p',U_s')^T$. Bounded solutions to Eq.~\ref{eq:wandereq1} therefore require the inhomogeneous terms to satisfy the Fredholm solvability condition, i.e., to be orthogonal to the nullspace of the adjoint operator $\mathcal{L}^*$, which is given by
\begin{equation}
    \mathcal{L}^*\begin{pmatrix}
        \Psi_e\\\Psi_p\\\Psi_s
    \end{pmatrix}=\begin{pmatrix}
        -\Psi_e+f'(U_e)[w_{ee}*\Psi_e+w_{pe}*\Psi_p+w_{se}*\Psi_s]\\-\Psi_p-f'(U_p)[w_{ep}*\Psi_e+w_{pp}*\Psi_p+w_{sp}*\Psi_s]\\-\Psi_s-f'(U_s)[w_{es}*\Psi_e+w_{ps}*\Psi_p+w_{ss}*\Psi_s]
    \end{pmatrix}.
\end{equation}
For the Heaviside firing-rate function, we find that $\ker(\mathcal{L}^*)$ is spanned by the adjoint null vector localized at the six threshold interfaces,
\begin{equation}\label{eq:wandereq2}
    \boldsymbol{\varphi}(x)=
    \begin{pmatrix}
        \delta(x+a_e)-\delta(x-a_e)\\
        \beta_p[\delta(x+a_p)-\delta(x-a_p)]\\
        \beta_s[\delta(x+a_s)-\delta(x-a_s)]
    \end{pmatrix},
\end{equation}
where $\beta_p$ and $\beta_s$ are given by
\begin{align*}
    \beta_p&=-\frac{[W_{ep}^--W_{ep}^+]+\beta_s[W_{sp}^--W_{sp}^+]}{|U_p'(a_p)|+W_{pp}^--W_{pp}^+},\\
    \beta_s&=\frac{[|U_p'(a_p)|+W_{pp}^--W_{pp}^+][|U_e'(a_e)|-W_{ee}^-+W_{ee}^+]+[W_{pe}^--W_{pe}^+][W_{ep}^--W_{ep}^+]}{[|U_p'(a_p)|+W_{pp}^--W_{pp}^+][W_{se}^--W_{se}^+]-[W_{pe}^--W_{pe}^+][W_{sp}^--W_{sp}^+]},
\end{align*}
with $W_{mn}^{\pm}=w_{mn}(a_m\pm a_n)$. Note that given the vector will span the nullspace, we have chosen to normalize the excitatory component to have unit amplitude to more easily solve for the relative $\beta_p$ and $\beta_s$.

Projecting the inhomogeneous terms in Eq.~~\ref{eq:wandereq1} onto the adjoint null vector $\boldsymbol{\varphi}$ and imposing the Fredholm solvability condition yields an effective stochastic equation for the bump displacement $d\Delta(t)$. Because the noise sources driving the three populations are independent and white in time, the projected increment has zero mean. Taking $\Delta(0)=0$, we therefore obtain $\langle \Delta(t)\rangle=0$, while the variance grows linearly in time as 
\begin{equation}
    \langle \Delta(t)^2 \rangle=\epsilon\frac{\mathscr{N}_e+\mathscr{N}_p+\mathscr{N}_s}{2[|U_e'(a_e)|+\tau_p \beta_p|U_p'(a_p)|+\tau_s \beta_s|U_s'(a_s)|]^2}t=\epsilon\mathscr{D}t
\end{equation}
with the components $\mathscr{N}_e=C_{ee}(0)-C_{ee}(2a_e)$, $\mathscr{N}_p=\beta_p^2[C_{pp}(0)-C_{pp}(2a_p)]$, and $\mathscr{N}_s=\beta_s^2[C_{ss}(0)-C_{ss}(2a_s)]$. Note that should one apply correlated noise across populations, then the cross correlation terms would also need to be included. 

\begin{figure}
    \centering
    \includegraphics[width=\linewidth]{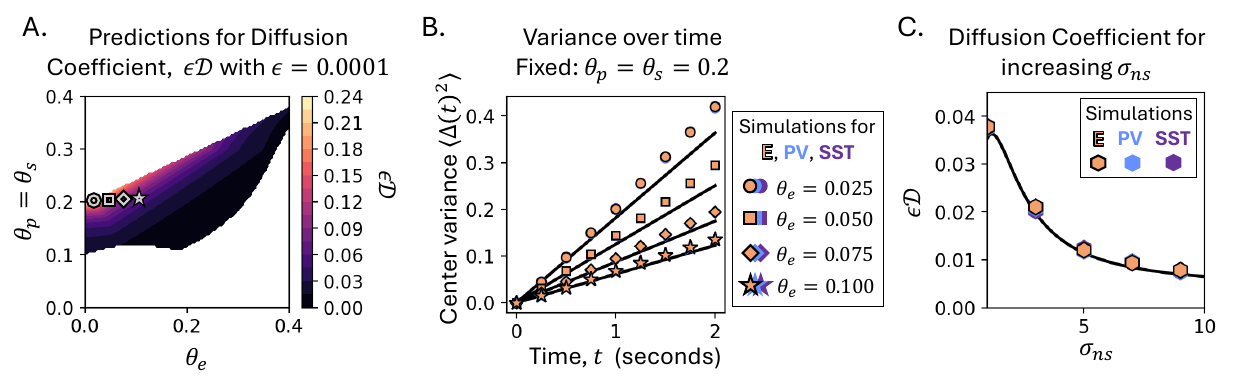}
    \caption{A. The scaled diffusion coefficient with respect to time in seconds and for stable solutions as shown in Figure \ref{fig:stabilityA}A. The four sampled thresholds were simulated for 2 seconds (200 steps in 10ms units) in panel B. Variance predictions (black lines) and variance over 10000 simulations each for increasing $\theta_e$ thresholds. While not always visible underneath the E population center variances, the variance of the centers for PV and SST were also shown. C. The theoretical scaled diffusion coefficient in units of degrees$^2$/second is shown in black for increasing synaptic footprint from the SST population (increasing long range connections). 10000 simulations were run for 1 seconds with the diffusion coefficient estimated from the slope between the final two sampled times. For A and B parameters are as in Table \ref{tab:model_params}. For C. $\sigma_{es}=\sigma_{ps}=\sigma_{ss}$ were increased simultaneously from 1 through 10 with all other parameters as in Figure \ref{fig:proportionStable}E-F.}
    \label{fig:diffusion}
\end{figure}

To test the diffusion approximation, we compare the theoretical prediction with direct simulations across a range of stable bump solutions (Fig.~\ref{fig:diffusion}). The predicted variance-growth rate $\epsilon\mathscr{D}$ varies substantially across threshold space (Fig.~\ref{fig:diffusion}A), and simulations at representative parameter values exhibit the expected linear growth in variance with close agreement to the theoretical slopes (Fig.~\ref{fig:diffusion}B). Thus, the reduced stochastic description accurately captures the long-time wandering of stable bumps in the weak-noise regime.

We now use the diffusion theory to examine how the spatial extent of SST inhibition influences stochastic memory dynamics. Increasing the spatial scales of projections from the SST population decreases the predicted diffusion rate $\epsilon\mathscr{D}$, with direct simulations showing the same trend (Fig.~\ref{fig:diffusion}C). Thus, broader SST connectivity not only alters the deterministic stability landscape described in Section~\ref{sec:weakIConnect}, but also reduces the rate at which stable memory representations wander under stochastic forcing. Together, these results demonstrate that inhibitory subtype structure can shape both the deterministic robustness and stochastic precision of persistent bump states.

\section{Discussion}
\label{sec:conclusion}
In this work, we developed a three-population neural field model to examine how distinct PV and SST inhibitory populations shape the dynamics of persistent activity bumps used to represent continuous memories. The Heaviside firing-rate formulation allowed us to characterize stationary bump solutions through their threshold interfaces and to reduce the corresponding linear stability problem into shifting and scaling subspaces. Using this framework, we found that population thresholds and inhibitory timescales can alter both the stability of bump solutions and the mode through which stability is lost, while changes in inhibitory connection strength and spatial scale can substantially reshape the stability landscape. In particular, broader inhibitory projections, consistent with the longer-range role assigned to SST neurons in our baseline architecture, expanded and deepened regions of stable bump activity. Weak connections among inhibitory populations could also have large, though nonmonotonic, effects on stability. Finally, projecting weak stochastic forcing onto the translational mode yielded an effective diffusion theory that accurately predicted bump wandering and showed that increasing the SST spatial footprint can reduce the diffusion of stable memory representations.

Several circuit mechanisms that distinguish cortical interneuron subtypes are not included in the present model. For example, we restricted the inhibitory population to PV and SST neurons, omitting additional subtypes such as vasoactive intestinal peptide (VIP) expressing  interneurons~\cite{campagnola2022local}. Incorporating VIP neurons would introduce another prominent disinhibitory pathway and an additional population with distinct spatial and temporal connectivity, potentially altering both the existence and stability of bump states. We also treated all synaptic connections as static, despite evidence that many of the connections considered here exhibit short-term facilitation or depression \cite{campagnola2022local, kullander2021cortical}. Activity-dependent synaptic efficacy could make the effective balance between PV and SST inhibition depend on the underlying network state, introducing additional timescales and potentially shifting the stability boundaries identified here. Extending the present interface-based analysis to include these mechanisms would provide a natural way to determine when the qualitative roles identified for inhibitory subtype structure persist in more dynamic circuit architectures.

Our analysis also focuses on a restricted class of memory states and perturbations. We considered stationary solutions in which each population possesses a single contiguous active region and examined stochastic fluctuations in the small-noise regime~\cite{li2021drifts}. This excludes multiple-bump solutions that could represent several simultaneously maintained items~\cite{krishnan2018synaptic}, as well as sufficiently large perturbations that move the system far from the neighborhood of a stationary bump~\cite{cihak2024robustly}. Such perturbations may arise from strong distractors or other transient inputs and could lead to bump displacement, annihilation, or transitions between distinct attractor configurations. In these regimes, the single-active-region construction and weak-noise diffusion reduction developed here need not remain valid. Characterizing the existence and interaction of multiple bump states, together with their responses to finite-amplitude perturbations, would therefore extend the present framework toward questions of memory competition and distractor-induced errors.

Another open question concerns population-specific wandering on intermediate timescales. Previous two-population neural field models have exhibited regimes in which excitatory and inhibitory bump positions develop different transient variances before converging to a common long-time diffusion rate \cite{cihak2022distinct}. The present stochastic reduction is constructed around the common translational zero mode and therefore describes the long-time diffusive motion shared by the three populations. Given our focus on the long-time dynamics and diffusion rate, we did not directly seek parameters for a comparable separation among the E, PV, and SST populations in the simulations, though we did observe this effect in some parameter sets. Transient relative displacements could nevertheless arise through the stable shifting modes identified in the linear stability analysis before those modes decay. A stochastic reduction that retains one or more of these stable modes alongside the translational mode could determine when population-specific wandering becomes appreciable and how its timescale depends on inhibitory subtype structure. 

The present results suggest that interneuron diversity is not simply an additional layer of biological detail in canonical excitatory-inhibitory bump models, but can play an important dynamical role in shaping persistent activity. At the same time, the resulting multi-population neural field provides a tractable testbed for dynamical-systems questions involving symmetry, bifurcations, interacting stability modes, and noise-driven motion along a neutral direction. Extending this framework to richer inhibitory circuits and state-dependent connectivity may therefore help identify both the mathematical mechanisms that organize these dynamics and the features of cortical inhibitory architecture most important for robust working-memory representations.

\newpage

\appendix
\section{Numerical Simulations}
\label{sec:numerics}
All code is written in Python 3. Key packages used were numpy for calculations, matplotlib for graphics, and scipy for some more calculations. The code used to generate the computational data and reproduce the figures of the paper will be deposited on Zenodo, an open research repository, upon paper’s acceptance.

Numerical simulations of Eqs.~\ref{eq: model} were carried out using forward Euler and Euler-Maruyama schemes for the deterministic system and stochastic systems, respectively. A single unit of time is considered to be 10ms and time increments of $dt=0.025$ or equivalently $0.25$ms were used to step through the system. The spatial axis $x$ is periodic and discretized with $n=2^{14}$ even increments. 

Convolutions on the ring were performed using Fast Fourier and Fast Inverse Fourier transforms from the scipy package. 
Active region halfwidths were found by using linear interpolation to find zero crossings over the bump profile minus the corresponding threshold. Eigenvalues of matrices were determined using the command linalg.eigvals(Matrix) from the numpy package.

The center of the bump was generally taken to be the peak, numerically identified via the numpy max function. In special cases a valid bump profile with a single connected active region would have a dip in the center, where the interfaces were then sampled and their midpoint was taken as the center of the bump.

\section*{Author Contributions}
Bilal Ahmed: Model formulation, Deterministic Calculations, Bifurcation Plots, Wandering Sample Plots, \& Code.
Heather Cihak: Problem Conceptualization, Model formulation, Figure Drafts and Finalization, Inhibitory Connection Exploration, Diffusion Calculation, Alternative population weightings, \& Code testing/review. Gregory Handy: Funding acquisition, Problem Conceptualization, Deterministic Calculations, Figure Drafts and Finalization, Calculation Review, \& Code testing/review. All authors contributed significantly to writing and editing of the manuscript, and have read and approved the final manuscript.
\bibliographystyle{siamplain}
\bibliography{references}

\end{document}